\documentclass{amsart}   % use for final version... turns off showkeys

\usepackage{amsmath,amssymb,latexsym, amscd}
\usepackage{exscale, cite, eps fig, graphics}

\renewcommand{\geq}{\geqslant}

\newcommand{\RM}{\mathbb{R}}
\newcommand{\ZM}{\mathbb{Z}}

\newcommand{\M}{\mathcal{M}}

\newtheorem{theorem}{Theorem}[section]

\newtheorem{proposition}[theorem]{Proposition}

\newtheorem*{main-theorem}{Main Theorem}
\newtheorem*{remark*}{Remark}
\numberwithin{equation}{section}

\title[The Whitham equation with capillarity and vorticity]
{Modulational instability in the Whitham equation with surface tension and vorticity}

\author[Hur]{Vera~Mikyoung~Hur}
\address{Department of Mathematics, University of Illinois at Urbana-Champaign, Urbana, IL 61801 USA}
\email{verahur@math.uiuc.edu}

\author[Johnson]{Mathew~A.~Johnson}
\address{Department of Mathematics, University of Kansas, Lawrence, KS 66045 USA} 
\email{matjohn@ku.edu}  

\date{\today}

\begin{document}

\maketitle

\begin{abstract}
We study modulational stability and instability in the Whitham equation, combining 
the dispersion relation of water waves and a nonlinearity of the shallow water equations, 
and modified to permit the effects of surface tension and constant vorticity.
When the surface tension coefficient is large, we show that 
a periodic traveling wave of sufficiently small amplitude is unstable to long wavelength perturbations
if the wave number is greater than a critical value, and stable otherwise, 
similarly to the Benjamin-Feir instability of gravity waves.
In the case of weak surface tension, we find intervals of stable and unstable wave numbers, 
whose boundaries are associated with the extremum of the group velocity, 
the resonance between the first and second harmonics, the resonance between long and short waves, 
and a resonance between dispersion and the nonlinearity. 
For each constant vorticity, we show that a periodic traveling wave of sufficiently small amplitude
is unstable if the wave number is greater than a critical value, and stable otherwise. 
Moreover it can be made stable for a sufficiently large vorticity. 
The results agree with those based upon
numerical computations or formal multiple-scale expansions to the physical problem. 
\end{abstract}

%%%%%%%%%%%%%%%%%%%%%%%%%%%%%%%%%%%%%%%%%%%%%%%%%
%%%%%%%%%%%%%%%%%%%%%%%%%%%%%%%%%%%%%%%%%%%%%%%%%
%%%%%%%%%%%%%%%%%%%%%%%%%%%%%%%%%%%%%%%%%%%%%%%%%
\section{Introduction}\label{sec:intro}
%%%%%%%%%%%%%%%%%%%%%%%%%%%%%%%%%%%%%%%%%%%%%%%%%
%%%%%%%%%%%%%%%%%%%%%%%%%%%%%%%%%%%%%%%%%%%%%%%%%
%%%%%%%%%%%%%%%%%%%%%%%%%%%%%%%%%%%%%%%%%%%%%%%%%

The Korteweg-de Vries (KdV) equation
\begin{equation}\label{E:KdV}
u_t+\sqrt{gd}\Big(1+\frac16d^2\partial_x^2\Big)u_x+uu_x=0
\end{equation}
approximates the water wave problem %(in the absence of the effects of surface tension and vorticity)
in a small amplitude and long wavelength regime,  
and furthermore, it satisfactorily explains solitary waves and cnoidal wave trains of the physical problem. 
Here, $t\in\mathbb{R}$ is proportional to elapsed time
and $x\in\mathbb{R}$ is the spatial variable in the predominant direction of wave propagation;
$u=u(x,t)$ is real-valued, related to the surface displacement from the undisturbed fluid depth $d$,
and $g$ is the constant due to gravitational acceleration. 
Throughout we express partial differentiation either by a subscript or using the symbol $\partial$. 
When waves are long compared to the fluid depth so that $kd\ll1$, $k$ the wave number,
one may expand the phase velocity\footnote{The dispersion relation of water waves is well-known, 
whose derivation dates back to the work of Airy~\cite{Airy} in 1845!}
%of a plane wave with the wave number $k$ 
for the water wave problem in the irrotational setting and write that
\begin{equation}\label{def:cW}
c_{WW}(k):=\pm\sqrt{\frac{g\tanh(kd)}{k}}
=\pm\sqrt{gd}\Big(1-\frac16k^2d^2\Big)+O(k^4d^4).
\end{equation}
Therefore the KdV equation may be regarded as to approximate up to second order
the dispersion relation of the water wave problem in the long wavelength regime.
Here and elsewhere, $\pm$ mean left and right moving waves.
(One may take, without loss of generality, the $+$ sign in \eqref{E:KdV}.
As a matter of fact, $x\mapsto -x$ and $u\mapsto -u$ turn one equation to the other.)

However the KdV equation misses peaking into sharp crests.
Furthermore waves in shallow water at times break into bores
whereas the KdV equation prevents singularity formation from solutions. 
This is not surprising since the phase velocity associated with the linear part of the KdV equation
poorly approximates\footnote{A relative error between $c_{WW}$ and the phase speed for the KdV equation
of, say, $10\%$ is made for $kd>1.242\dots$.} $c_{WW}$ when $kd$ becomes large.
One may incidentally argue that peaking and breaking are high frequency phenomena, 
for which short wave components become important 
and the long wavelength assumption is no longer adequate.

Whitham therefore emphasized in \cite{Whitham} that 
``it is intriguing to know what kind of simpler mathematical equation (than the physical problem)
could include" breaking and peaking, and he put forward 
\begin{equation}\label{E:whitham}
u_t+\mathcal{M}u_x+uu_x=0,
\end{equation}
where $\mathcal{M}$ is a Fourier multiplier, defined via its symbol as
\begin{equation}\label{def:m}
\widehat{\mathcal{M}f}(k)=m(k)\widehat{f}(k),
\end{equation}
and $m=c_{WW}$ (see \eqref{def:cW}).
It combines the full range of dispersion of water waves, rather than a second order approximation, 
and the nonlinearity of the shallow water equations, 
and hence it may offer an improvement over the KdV equation for short and intermediately long waves.
As a matter of fact, numerical experiments in \cite{BKN} and \cite{MKD} indicate that 
Whitham's model approximates waves in water on par with or better 
than the KdV or other shallow water equations do in some respects outside the long wavelength regime. 
(The KdV equation, on the contrary, seems a better approximation in the long wavelength regime.)
%As a matter of fact, numerical experiments in \cite{BKN} and \cite{MKD} 
%indicate that Whitham's model approximates waves in water on par with or better
%than the KdV and other shallow water equations do in a wide range of amplitude and wavelength parameters. 
%END attempt
Moreover Whitham conjectured that \eqref{E:whitham}-\eqref{def:m}, where $m=c_{WW}$,
would explain breaking and peaking. 
Wave breaking --- bounded solutions with unbounded derivatives --- 
in the Whitham equation for water waves was recently settled in \cite{Hur-breaking1}.
(See \cite{NS} and \cite{CE-breaking} for related results, and the discussion in \cite{Hur-breaking1}.)
Sharp crests in its periodic traveling wave of maximum amplitude
was numerically supported in \cite{EK2} and an analytical result was announced in \cite{Mats}.
(A complete proof of global bifurcation will be reported elsewhere.) 
In the last few years, the Whitham equation for water waves 
received renewed attention from a number of other vantage points; 
see \cite{EK}, \cite{EGW}, \cite{HJ2}, \cite{SKCK2014}, \cite{EEP}, among others.

\

Benjamin and Feir in \cite{BF} and, independently, Whitham in \cite{Whitham1967},
on the other hand, formally argued that 
a $2\pi/k$-periodic traveling wave in water under the influence of gravity would be unstable, 
leading to sidebands growth, namely the Benjamin-Feir or modulational instability,
provided that $kd>1.363\dots$.
Benny and Newell in \cite{BN} and Zakharov in \cite{Zakharov-WW},
and later Hasimoto and Ono in \cite{HO} formally derived a nonlinear Schr\"odinger equation, 
which describes the long time behavior of the envelope of Stokes waves,
and deduced the Benjamin-Feir instability. 
Kawahara in \cite{kawahara}, Djordjevi\'c and Redekopp in \cite{DR} and later Hogan in \cite{Hogan}
extended the results to capillary-gravity waves. 

The Benjamin-Feir instability is another high frequency phenomenon in water waves, 
which does not manifest in the KdV and other shallow water equations\footnote{
Recently in \cite{HP1}, a small-amplitude, $2\pi/k$-periodic traveling wave 
of the Benjamin-Bona-Mahony equation was shown to be 
modulationally unstable if $k>\sqrt{3}$, and stable otherwise. But the instability mechanism is different
from that in the Whitham equation for water waves or the water wave problem.}. 
Including all the dispersion present in water waves, notwithstanding, 
the Whitham equation may predict short waves' instabilities. 
Recently in \cite{HJ2}, as a matter of fact,  analytically proved was that 
a small-amplitude, $2\pi/k$-periodic traveling wave of \eqref{E:whitham}-\eqref{def:m}, where $m=c_{WW}$,
be unstable to long wavelength perturbations if $kd>1.145\dots$, 
and spectrally stable to square integrable perturbations otherwise. 
Numerical experiments in \cite{SKCK2014} bear it out.  
Note that the critical wave number compares reasonably well with that in \cite{Whitham1967}. 

\

Here we take matters further and examine the effects of surface tension and constant vorticity 
upon the modulational stability and instability in the Whitham equation for water waves. 
In the presence of effects of surface tension, we extend the results in \cite{HJ2} 
to \eqref{E:whitham}-\eqref{def:m}, where
\begin{equation}\label{def:mT}
m(k)=\pm\sqrt{(g+Tk^2)\frac{\tanh(kd)}{k}}
\end{equation}
is the phase velocity for capillary-gravity waves in the irrotational setting; 
$T\geq 0$ is the coefficient of surface tension.
At the air-sea interface, $g\approx9.81m\cdot s^{-2}$ and $T\approx7.3\times10^{-3}N\cdot m^{-1}$,
and thus the effects of surface tension become negligible 
if wavelengths are several times greater than $2\pi\sqrt{T/g}\approx1.7cm$.
In the case of $T=0$, \eqref{def:mT} reduces to \eqref{def:cW}.
 
We show that the effects of surface tension considerably alters the modulational stability and instability 
of small-amplitude, periodic traveling waves of \eqref{E:whitham}-\eqref{def:m} and \eqref{def:mT}.
In the $kd$ versus $k\sqrt{T/g}$ plane, specifically,
we find several regions of modulational stability and instability,
whose boundaries are associated with the extremum of the group velocity, 
the resonance between the fundamental mode and the second harmonic,
the resonance between the limiting long wave and a short wave, 
and a resonance between dispersion and the nonlinearity, 
which in the absence of the effects of surface tension takes place at $kd=1.145\dots$. 
See Figure~\ref{f:contour}. 
Our results by and large reproduce those in \cite{K} and \cite{DR}, for instance,
via formal multiple-scale expansions to the capillary-gravity wave problem.
In contrast, here we make a rigorous calculation of long wavelengths perturbations
to a formal approximate model of the problem.
These two approaches compare (fortuitously) reasonably well. 

When the coefficient of surface tension is large, or equivalently, the fluid depth is small so that $T/gd^2>1/3$, 
we show that a small-amplitude, periodic traveling wave 
of \eqref{E:whitham}-\eqref{def:m} and \eqref{def:mT} is unstable to long wavelength perturbations, 
provided that the wave number is greater than a critical value, 
and stable otherwise, similarly to the gravity wave setting; see \cite{HJ2}. 
In the case of $0<T/gd^2<1/3$, on the contrary, 
we find three intervals of modulationally stable wave numbers and three intervals of unstable wave numbers.

In the capillary wave setting, i.e. $g=0$, furthermore, we deduce that 
a small-amplitude, periodic traveling wave of \eqref{E:whitham}-\eqref{def:m}, \eqref{def:mT} 
is modulationally unstable regardlessly of the wave number, 
which is consistent with the results in \cite{kawahara} and \cite{DR}, for instance.
In the deep water limit as $kd\to\infty$, we deduce that 
it is modulationally stable if $2/\sqrt{3}-1<Tk^2/g<1/2$, and unstable otherwise,
which is consistent with the results in \cite{Lighthill-deep} and \cite{Simmons}, for instance,
via formal multiple-scale expansions to the physical problem.

\

In the presence of the effects of constant vorticity, similarly,
we extend the results in \cite{HJ2} to \eqref{E:whitham}-\eqref{def:m}, where
\begin{equation}\label{def:mw}
m(k)=\frac{\gamma\tanh(kd)}{2k}\pm\sqrt{\frac{g\tanh(kd)}{k}+\frac{\gamma^2\tanh^2(kd)}{4k^2}}
\end{equation}
is the phase velocity of Stokes waves with the vorticity $\gamma\in\mathbb{R}$ a constant;
see\footnote{This seminal work quickly led to a flurry of research activities about 
traveling water waves with vorticity, which we do not do justice here.}
 \cite{CS2004}, for instance. %\cite{Wah09}, \cite{CV11}, and \cite{CSV15}.
For general non-constant vorticities, the dispersion relation may not be written in closed form;
see \cite{HL08} and \cite{Kar12}, for instance.

The zero vorticity setting serves as a reasonable approximation in some circumstances
and it is easier to handle, analytically and numerically,  than the rotational counterpart.
But the inclusion of vorticities is likely to be far more applicable. 
As a matter of fact, waves are most commonly generated by wind 
and in any region where the wind is blowing there is a surface drift of the water. 
Moreover currents cause shear at the bed of the sea or of a river. 
Constant vorticity is representative 
when waves are short compared to the length scale of the vorticity distribution; see \cite{SP}, for instance.
%For example, tidal currents are well modeled in \cite{Con11} by constant vorticities.  For discussions regarding the physical relevance of flows with constant vorticity, see \cite{SP88} and \cite{AJ15}.
%Of course, constant vorticities are perhaps a poor choice for wind-driven currents and other flows. 

Colin, Dias and Ghidaglia in \cite{CDG95} formally derived an amplitude equation 
for a modulated Stokes wave with vorticity, not necessarily constant, 
which is the same as that in the irrotational setting, 
whereby they argued that vorticity does not qualitatively alter 
the modulational stability and instability of Stokes waves.
(The velocity field depends upon the vorticity, though.) 
But vorticity does influence upon the stability of water waves. 
As a matter of fact, some small-amplitude, periodic traveling waves over a shear flow 
were shown in \cite{HL08} to be linearly unstable to the same period perturbations, 
regardlessly of the vorticity strength. 
Incidentally constant vorticity does not seem to invite instability of the kind.

For each $\gamma\in\mathbb{R}$ we show that 
a small-amplitude, periodic traveling wave of \eqref{E:whitham}-\eqref{def:m} and \eqref{def:mw}
is unstable to long wavelength perturbations, provided that the wave number is greater than a critical value, 
and stable otherwise, similar to the zero vorticity setting; see \cite{HJ2}. 
But the critical wave number depends upon the sign and the strength of the constant vorticity.
Specifically, the critical wave number increases unboundedly 
as $\gamma\to-\infty$ in the case of the $+$ sign in \eqref{def:mw} 
and as $\gamma\to+\infty$ in the case of the $-$ sign; see Figure~\ref{f:vorcontour}. 
For each $k>0$, therefore, 
a small-amplitude, $2\pi/k$-periodic traveling wave of \eqref{E:whitham}-\eqref{def:m}, \eqref{def:mw}
is modulationally stable, for constant vorticities sufficiently large in strength,
whose sign depends upon the sign in \eqref{def:mw}.
Our result helps to shed light into the stabilization effects of constant vorticities upon Stokes waves,
which recent numerical experiments in \cite{OSV}, for instance, suggest.

%Recent numerical experiments in \cite{OSV} suggest that a small-amplitude Stokes wave in water can be made modulationally stable with a sufficiently large constant vorticity. 

\subsection*{Notation}
Let $L^p_{2\pi}$ in the range $p\in[1,\infty]$ denote the space of 
$2\pi$-periodic, measurable, real or complex valued functions over $\RM$ such that
\[
\|f\|_{L^p_{2\pi}}=\left(\frac{1}{2\pi}\int_{-\pi}^\pi |f|^p~dx\right)^{1/p}<+\infty\quad\textrm{if}\quad p<\infty
\]
and essentially bounded if $p=\infty$. 
Let $H^1_{2\pi}$ denote the space of $L^2_{2\pi}$-functions such that $f'\in L^2_{2\pi}$,
and $H^\infty_{2\pi}=\cap_{k=0}^\infty H^k_{2\pi}$.

We express $f\in L^1_{2\pi}$ as a Fourier series:
\[
f(z)\sim\sum_{n\in\ZM}\widehat{f}_n~e^{inz},\quad\textrm{where}\quad 
\widehat{f}_n=\frac{1}{2\pi}\int_{-\pi}^\pi e^{-inz}f(z)dz.
\]
If $f\in L^p_{2\pi}$, $p>1$, then the Fourier series converges to $f$ pointwise almost everywhere.
We define the half-integer Sobolev space $H^{1/2}_{2\pi}$ via the norm
\[
\|f\|_{H^{1/2}_{2\pi}}^2=\widehat{f}_0^2+\sum_{n\in\ZM}|n||\widehat{f}_n|^2.
\]

%%%%%%%%%%%%%%%%%%%%%%%%%%%%%%%%%%%%%%%%%%%%%%%%%
%%%%%%%%%%%%%%%%%%%%%%%%%%%%%%%%%%%%%%%%%%%%%%%%%
%%%%%%%%%%%%%%%%%%%%%%%%%%%%%%%%%%%%%%%%%%%%%%%%%
\section{Small amplitude capillary-gravity waves}\label{sec:existence}
%%%%%%%%%%%%%%%%%%%%%%%%%%%%%%%%%%%%%%%%%%%%%%%%%
%%%%%%%%%%%%%%%%%%%%%%%%%%%%%%%%%%%%%%%%%%%%%%%%%
%%%%%%%%%%%%%%%%%%%%%%%%%%%%%%%%%%%%%%%%%%%%%%%%%
We begin by discussing how one constructs periodic traveling waves 
of the {\em capillarity-modified Whitham equation}, \eqref{E:whitham}-\eqref{def:m} and, abusing notation, 
\begin{equation}\label{defT:m}
m(k)=\sqrt{gd\Big(1+\frac{T}{gd^2}k^2\Big)~\frac{\tanh k}{k}}.
\end{equation}
Moreover we record their small amplitude asymptotics.
Note that %As a matter of fact, 
$kd \mapsto k$ turns \eqref{def:mT} to \eqref{defT:m}. 
One may assume, without loss of generality, the $+$ sign in \eqref{def:mT}. 
As a matter of fact, $x\mapsto -x$ and $u\mapsto -u$ turn one equation to the other. 
%Moreover, below we record their small amplitude asymptotics. 

In the long wave limit as $k\to 0$, one may expand the symbol in \eqref{defT:m} and write that
\[
m(k)=\sqrt{gd}\Big(1-\frac12\Big(\frac13-\frac{T}{gd^2}\Big)k^2\Big)+O(k^4).
\]
Unless $T/gd^2=1/3$, therefore, the KdV equation with surface tension
\begin{equation}\label{ET:KdV}
u_t+\sqrt{gd}\Big(1+\frac12\Big(\frac13-\frac{T}{gd^2}\Big)\partial_x^2\Big)u+uu_x=0
\end{equation}
approximates up to second order 
the dispersion of the capillarity-modified Whitham equation, 
and hence the water wave problem with surface tension.
%As a matter of fact, solutions of \eqref{ET:KdV} and \eqref{E:whitham}-\eqref{def:m}, \eqref{defT:m}
%exist and they converge to those of the water wave problem with surface tension at the order of $O(k^2)$
%during a relevant interval of time. 
In the short wave limit $k\to \infty$, on the other hand, 
$m(k)\sim \sqrt{T|k|}$, and \eqref{E:whitham}-\eqref{def:m} and \eqref{defT:m} 
become, after normalization of parameters, a perturbation of 
\[
u_t+\Lambda^{1/2}u_x+uu_x=0, \qquad \text{where}\quad \widehat{\Lambda f}(k)=|k|\widehat{f}(k),
\]
in the family of the KdV equations with fractional dispersion, which incidentally is $L^2$-critical.

\

A traveling wave solution of \eqref{E:whitham} takes the form $u(x,t)=u(x-ct)$, 
where $c>0$ and $u$ satisfies by quadrature that 
\begin{equation}\label{quad}
\mathcal{M} u-cu+u^2=(m(0)-c)^2b
\end{equation}
for some $b\in\RM$. In other words, it propagates at a constant speed $c$ without change of shape.
We seek solutions of \eqref{quad} of the form
\[
u(x)=w(z),\qquad z=k x,
\]
where $k>0$ is the wave number and $w$ is $2\pi$-periodic, satisfying that 
\begin{equation}\label{quad2}
\mathcal{M}_{k}w-cw+w^2=(m(0)-c)^2b.
\end{equation}
Here and elsewhere,
\[
\mathcal{M}_{k}e^{inz}=m(k n)e^{inz}\quad\textrm{for}\quad n\in\mathbb{Z}
\]
and it is extended by linearity and continuity. 

For each $T>0$, note that $m$ is smooth, even, 
strictly and unboundedly increasing away from an interval of low wave numbers.
For each $T>0$ and $k>0$, therefore, $\mathcal{M}_k$ maps even functions to even functions.
%For each $T>0$, moreover, 
%Moreover, for each $T>0$ the operator
Moreover
\[
\mathcal{M}_{k}:H^s_{2\pi}\to H^{s+1/2}_{2\pi}\quad\textrm{for all $k>0$}\quad\text{for all $s\geq 0$}
\]
is bounded. Therefore, if $w\in H^1_{2\pi}$ satisfies \eqref{quad2} 
for some $c>0$, $k>0$ and $b\in\mathbb{R}$ then $w\in H^\infty_{2\pi}$.
Indeed $(\mathcal{M}_k-c)w=-w^2+(m(0)-c)^2b \in H^1_{2\pi}$ by the Sobolev inequality,
and hence $w\in H^{3/2}_{2\pi}$. The claim then follows from a bootstrap argument. 
Thanks to Galilean invariance of \eqref{quad2} under
\begin{equation}\label{E:galilean}
w(z)\mapsto w(z)+v\quad\text{and}\quad c\mapsto c-2v, \quad
(m(0)-c)^2b\mapsto (m(0)-c)^2b-cv+v^2
\end{equation}
for any $v\in\mathbb{R}$, furthermore, we may take, without loss of generality, $b=0$ in \eqref{quad2}.
In the following section, however, we shall make use of variations of the carrier wave in the $b$-variable. 

\

For each $T>0$ and $k>0$ fixed, clearly, 
the zero function is a solution of \eqref{quad2}, where $b=0$, for all $c>0$.
%Fixing $b=0$, clearly the zero function is a solution 
%of \eqref{quad2} for all $c>0$ for each fixed $T>0$ and $k>0$. %, where $b=0$, for all $c>0$. 
We shall determine at which values of $c$ a branch of non-constant even $H^1_{2\pi}$-solutions, 
and hence smooth solutions, emanates from the zero solution. 
It follows from bifurcation theory that a necessary condition is that 
the associated linearized operator admits a nontrivial kernel. 
While this is not in general a sufficient condition, bifurcaiton indeed occurs
if the kernel is one dimensional.
%This is not in general a sufficient condition, 
%but bifurcation indeed takes place if the kernel is one dimensional. 

%NEW-MJ  Not quite... could have more in kernel... reworded to concentrate on primary branch only
%A straightforward calculation reveals that
%\[
%\ker(\mathcal{M}_k-c)=\text{span}\{\cos(nz)\}
%\]
%in the sector of even functions in $H^1_{2\pi}$, 
%provided that $c=m(kn)$, $n\in\mathbb{N}$, and the kernel is trivial otherwise. 
A straightforward calculation reveals that the kernel of the linear part $\mathcal{M}_k-c$ of \eqref{quad2}
is trivial unless $c=m(kn)$ for some $n$ an integer, in which case 
\[
\cos(n z)\in\ker(\mathcal{M}_k-c)
\]
in the sector of even functions in $H^1_{2\pi}$. In the case of $n=1$, for instance, 
bifurcation is guaranteed to occur, provided that $m(k)\neq m(kN)$ for all $N=2, 3, \dots$.  
Otherwise, $m(k)=m(kN)$ for some $N\geq 2$ an integer, and
the kernel at $c=m(k)$ is two dimensional, generated by $\cos(z)$ and $\cos(Nz)$. 
Physically, it means that the phase velocities of the fundamental mode and a higher harmonic coincide.
In what follows, we take $n=1$ and examine when bifurcation occurs at $c=m(k)$.
%We now consider when a bifurcation is guaranteed to occur from $c=m(k)$.

\begin{figure}[h]
(a)~~\includegraphics[scale=0.32]{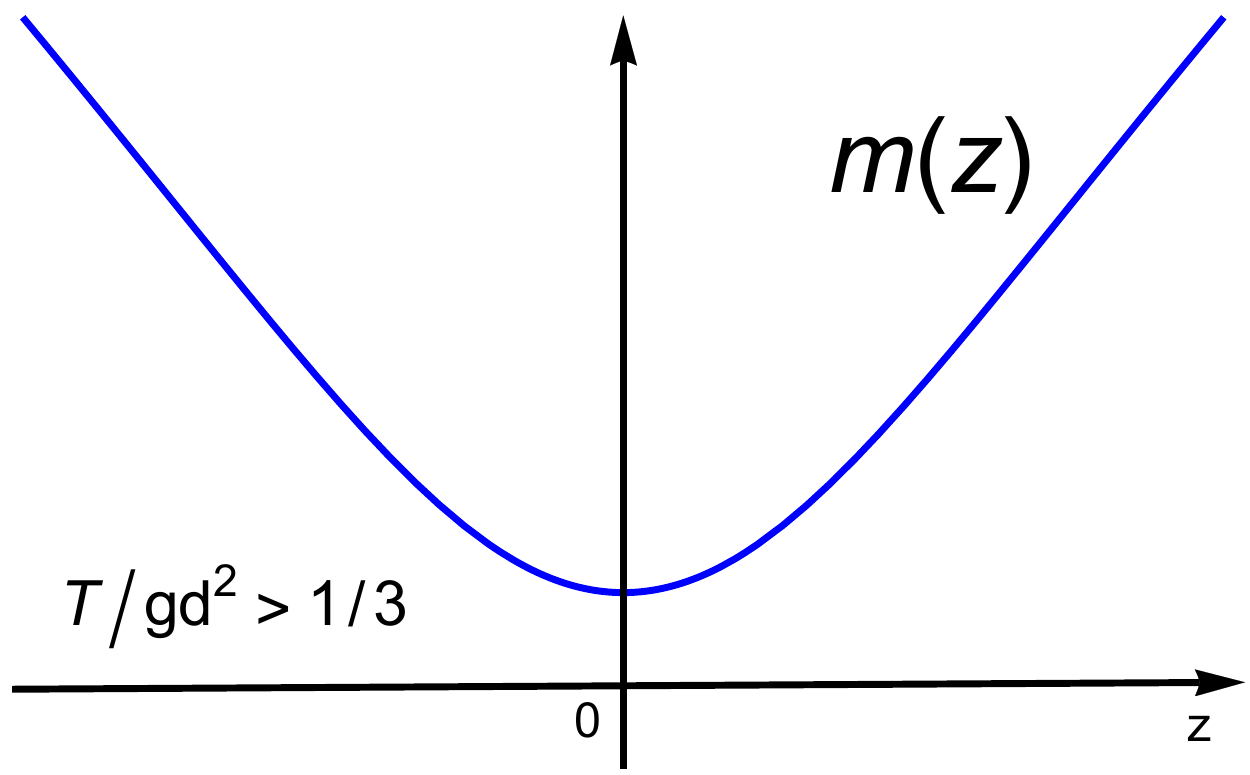}\quad
(b)~~\includegraphics[scale=0.32]{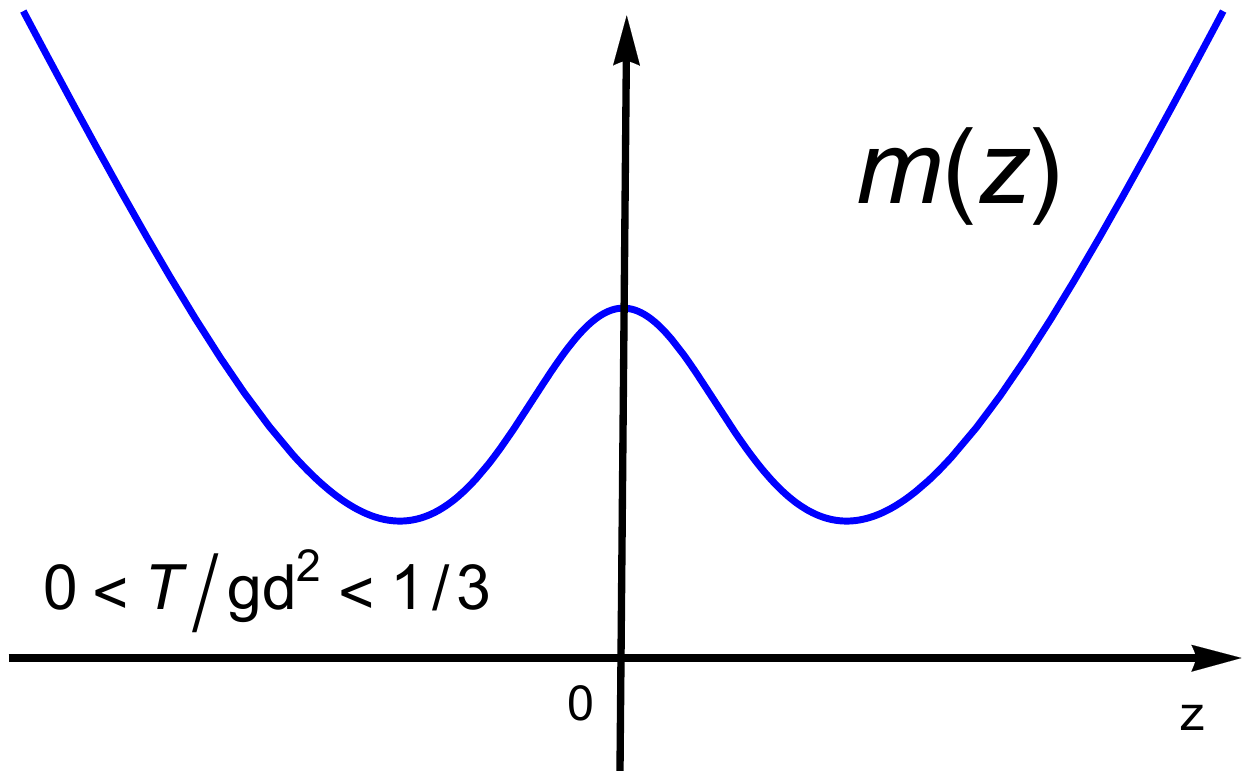}\quad
(c)~~\includegraphics[scale=0.32]{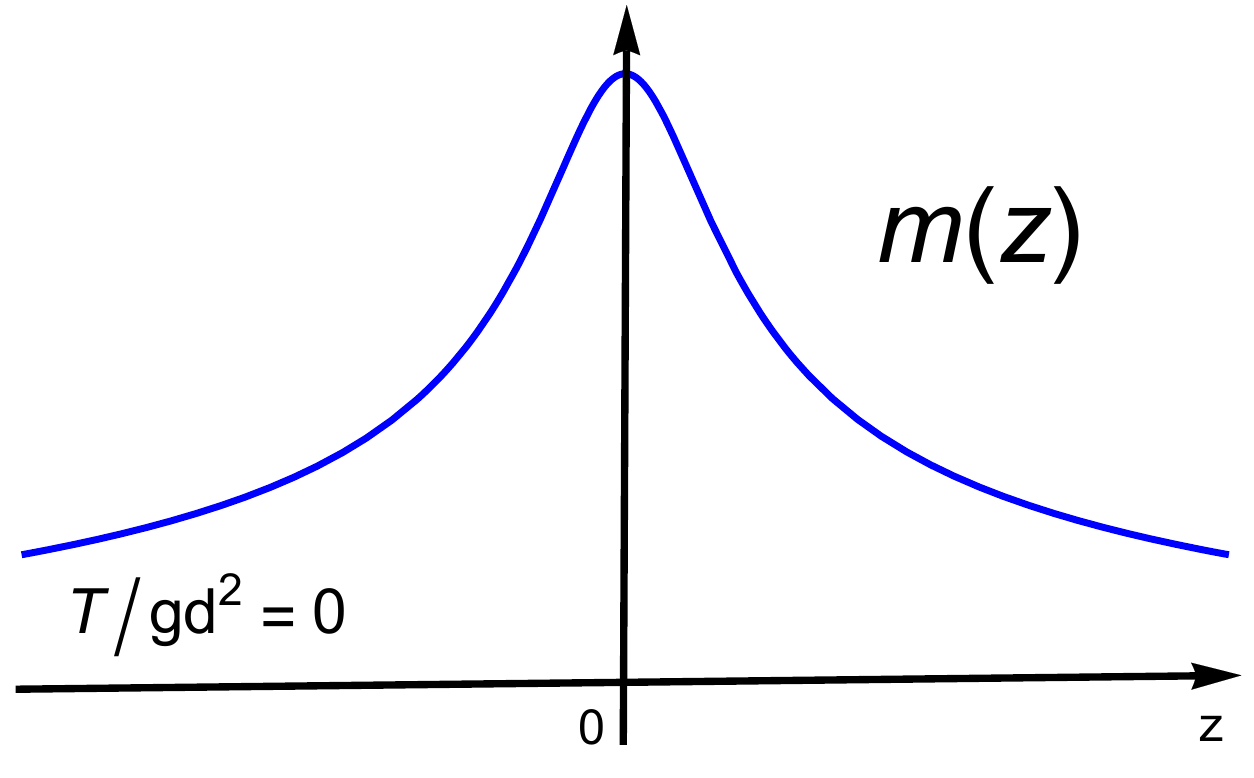}
\caption{Schematic plots of the symbol $m$ in \eqref{defT:m} for (a) $T/gd^2>1/3$, (b) $0<T/gd^2<1/3$, and (c) $T=0$.}\label{f:capsymbol}
\end{figure}

When the coefficient of surface tension is large, or equivalently,
the fluid depth is small so that $T/gd^2\geq1/3$,  it is readily verifed that
$m(z)>0$ for $z>0$ and it monotonically increases, unbounded above; see Figure~\ref{f:capsymbol}(a). 
Therefore $m(kN)>m(k)$ for all $k>0$ for all $N=2,3,\dots$. 
In other words, bifurcation indeed occurs from the zero solution at $c=m(k)$.%each $c=m(kn)$ .

%%%%%%%%%%%%%%%%%%%%%%%%%%%%%%%%%%%%%%%%%%%%%%%%%
%%%%%%%%%%%%%%%%%%%%%%%%%%%%%%%%%%%%%%%%%%%%%%%%%
%%%%%%%%%%%%%%%%%%%%%%%%%%%%%%%%%%%%%%%%%%%%%%%%%
%%%%%%%%%%%%%%%%%%%%%%%%%%%%%%%%%%%%%%%%%%%%%%%%%

In the case of $0<T/gd^2<1/3$, on the contrary, 
$m'(0)=0$, $m''(0)<0$ and $m(z)\to\infty$ as $z\to \infty$. 
In particular, the graph of $m$ supports a single local minimum over $(0,\infty)$; 
see Figure~\ref{f:capsymbol}(b). 
A straightforward calculation then reveals that 
for each $N\geq 2$ an integer, one finds a unique $k_N>0$ such that $m(k_N)=m(k_NN)$.
Moreover $k_N$ strictly decreases to zero as $N\to \infty$; see \cite[Lemma~3]{M2013}, for instance.
In the case of $T=0$, furthermore, $m(z)>0$ for $z>0$, and it monotonically decreases to zero as $z\to\infty$;
see Figure~\ref{f:capsymbol}(c). Therefore $m(k)>m(kN)$ for all $k>0$ for all $N=2,3,\dots$ 
and bifurcation occurs at %each $c=m(kn)$; 
$c=m(k)$; see \cite{HJ2}, for instance.

Let 
\begin{equation}\label{def:sigma}
\Sigma=\begin{cases}
(0,\infty)\quad &\text{if}\quad T/gd^2\geq1/3, \\
(0,\infty)\setminus \{k_N:N=2,3,\dots\}\quad&\text{if}\quad 0<T/gd^2<1/3.
\end{cases}
\end{equation}
For each $k\in\Sigma$ one may employ bifurcation theory
and obtain a one-parameter family of smooth and even solutions of \eqref{quad2} 
near the zero solution and $c=m(k)$.
We summarize the conclusion below. The detail is found in \cite{HJ2}, for instance.
We merely pause to remark that for $|b|$ sufficiently small
we appeal to the Galilean invariance in \eqref{E:galilean}.

%%%%%%%%%%%%%%%%%%%%%%%%%%%%%%%%%%%%%%%%%%%%%%%%%
%%%%%%%%%%%%%%%%%%%%%%%%%%%%%%%%%%%%%%%%%%%%%%%%%
\begin{proposition}[Existence]\label{p:existence} Let $T>0$. 
For each $k\in \Sigma$ and for $|b|$ sufficiently small, 
a family of periodic traveling waves of \eqref{E:whitham}-\eqref{def:m}, \eqref{defT:m} exists and
\[
u(x,t)=w(a,b)(k(x-c(k,a,b)t))=:w(k,a,b)(z)
\]
for $|a|$ sufficiently small, where $w$ and $c$ depend smoothly upon $k$, $a$, $b$. 
Moreover $w$ is smooth, even and $2\pi$-periodic in $z$, and $c$ is even in $a$.  Furthermore, 
\begin{align}
w(k,a,b)(z)=&w_0(k,b)+a\cos(z) \notag \\
&+\frac{1}{2}a^2\Big(\frac{1}{m(k)-m(0)}+\frac{\cos(2z)}{m(k)-m(2k)}\Big)+O(a(a^2+b^2))\label{E:w}
\intertext{and}
c(k,a,b)=&c_0(k,b)+a^2\Big(\frac{1}{m(k)-m(0)}+\frac{1}{2}\frac{1}{m(k)-m(2k)}\Big)+O(a(a^2+b^2))\label{E:c}
\end{align}
as $|a|,|b|\to 0$, where
\begin{align*}
c_0(k,b)=&m(k)+2(m(0)-m(k))(b-3b^2+O(b^3))
\intertext{and}
w_0(k,b)=&(m(0)-m(k))(b-b^2+O(b^3)).
\end{align*}
\end{proposition}
%%%%%%%%%%%%%%%%%%%%%%%%%%%%%%%%%%%%%%%%%%%%%%%%%
%%%%%%%%%%%%%%%%%%%%%%%%%%%%%%%%%%%%%%%%%%%%%%%%%

In case $k\notin\Sigma$, i.e. $m(k)=m(kN)$ for some $n\geq N$ an integer, 
the classical theory of bifurcation does not apply, but one may nevertheless employ 
the Lyapunov-Schmidt method, for instance, and obtain a branch of small amplitude solutions.
See \cite{Jones}% and \cite{MM13}
, for instance, for the detail in the capillary-gravity wave problem. 

In particular, Wilton ripples arise 
when the phase velocities of the fundamental mode and the second harmonic match, i.e. 
$m(k)=m(2k)$, or equivalently, 
\[
Tk^2/g=\frac{\tanh^2(kd)}{3-\tanh^2(kd)},
\] 
which in the deep water limit as $kd\to \infty$ becomes $Tk^2/g=1/2$. 
The stability theory in the following section, or in \cite{HJ2}, does not apply for resonant solutions, 
however, and it is an interesting direction of future investigations.

%%%%%%%%%%%%%%%%%%%%%%%%%%%%%%%%%%%%%%%%%%%%%%%%%
%%%%%%%%%%%%%%%%%%%%%%%%%%%%%%%%%%%%%%%%%%%%%%%%%
%%%%%%%%%%%%%%%%%%%%%%%%%%%%%%%%%%%%%%%%%%%%%%%%%
\section{Modulational Instability Index}\label{sec:MI}
%%%%%%%%%%%%%%%%%%%%%%%%%%%%%%%%%%%%%%%%%%%%%%%%%
%%%%%%%%%%%%%%%%%%%%%%%%%%%%%%%%%%%%%%%%%%%%%%%%%
%%%%%%%%%%%%%%%%%%%%%%%%%%%%%%%%%%%%%%%%%%%%%%%%%

Let $T>0$. Let $w=w(k,a,b)$ and $c=c(k,a,b)$, 
for $k\in\Sigma$ and $|a|$, $|b|$ sufficiently small, form 
a small-amplitude, $2\pi/k$-periodic traveling wave of \eqref{E:whitham}-\eqref{def:m} and \eqref{defT:m},
whose existence follows from the previous section.
We adapt the developments in \cite{HJ2} (in the gravity wave setting, i.e. $T=0$) 
and study its modulational stability and instability. 

\

Linearizing \eqref{E:whitham} about $w$ in the frame of reference moving at the speed $c$, we arrive at that
\[
v_t+k\partial_z(\mathcal{M}_{k}-c+2w)v=0.
\]
Seeking a solution of the form $v(z,t)=e^{\lambda k t}v(z)$, 
$\lambda\in\mathbb{C}$ and $v\in L^2(\mathbb{R})$, 
moreover, we arrive at the spectral problem
\begin{equation}\label{spec}
\lambda v=\partial_z(-\mathcal{M}_{k}+c-2w)v=:\mathcal{L}(k,a,b)v.
\end{equation}
We say that $w$ is \emph{spectrally unstable} if
the $L^2(\RM)$-spectrum of $\mathcal{L}$ intersects the open, right half plane of $\mathbb{C}$
and it is (spectrally) stable otherwise. Note that $v$ needs not have the same period as $w$.
Since the spectrum of $\mathcal{L}$ is symmetric
with respect to reflections about the real and imaginary axes,
$w$ is spectrally unstable if and only if 
the $L^2(\mathbb{R})$-spectrum of $\mathcal{L}$ is {\emph not} contained in the imaginary axis.

It follows from Floquet theory (see \cite{BHJ}, for instance, and references therein) that 
nontrivial solutions of \eqref{spec} cannot be integrable over $\mathbb{R}$.
They are at best bounded over $\mathbb{R}$,
and the $L^2(\mathbb{R})$-spectrum of $\mathcal{L}$ is purely essential.
In the case of, e.g., the KdV equation, the essential spectrum of the associated linearized operator
may be related to eigenvalues of a monodromy map,
and the stability may be studied with the help of Evans function techniques and other ODE methods.
Confronted with {\em nonlocal} operators, ODE methods may not be readily applicable.
It follows from Floquet theory (see \cite{BHJ}, for instance, and references therein), instead, 
that any eigenfunction of \eqref{spec} takes the form 
\[
v(z)=e^{i\xi z}\phi(z),\qquad\text{where $\phi$ is $2\pi$-periodic and $\xi\in[-1/2,1/2)$.}
\]
Therefore $\lambda\in\mathbb{C}$ belongs to the $L^2(\mathbb{R})$-spectrum of $\mathcal{L}$
if and only if 
\[
\lambda v=\mathcal{L}v\quad\text{and}\quad v(z+2\pi)=e^{2\pi i\xi}v(z),
\]
or equivalently,
\begin{equation}\label{E:Lxi}
\lambda\phi=e^{-i\xi z}\mathcal{L}(k,a,b)e^{i\xi z}\phi=:\mathcal{L}_\xi(k,a,b)\phi
\end{equation}
for some $\xi\in[-1/2,1/2)$ for some $2\pi$-periodic $\phi$. 
For each $\xi\in[-1/2,1/2)$ the $L^2_{2\pi}$-spectrum of $\mathcal{L}_\xi$ 
comprises entirely of discrete eigenvalues of finite multiplicities. Furthermore
\[
\text{spec}_{L^2(\RM)}(\mathcal{L}(k,a,b))=
\bigcup_{\xi\in[-1/2,1/2)}\text{spec}_{L^2_{2\pi}}(\mathcal{L}_\xi(k,a,b)).
\]
Thereby we parametrize the essential $L^2(\mathbb{R})$-spectrum of $\mathcal{L}$
by the one-parameter family of point $L^2_{2\pi}$-spectra of $\mathcal{L}_\xi$'s.

Note that $\xi=0$ corresponds to $2\pi$-periodic perturbations, 
i.e. perturbations with the same period as the carrier wave,
and $|\xi|\ll 1$ corresponds to long wavelength perturbations %or slow modulations 
of the carrier wave.  
%NEW-MJ
Slow modulations of the carrier wave are a special class of long wavelength perturbations 
in which the effects of the perturbation are to slowly vary, namely modulate, the wave characteristics 
--- $k$, $a$, $b$, in the present setting --- and the spatial translation mode.  
Note that variations in these parameters provide spectral information 
about $\mathcal{L}$ at the origin in the spectral plane; see \cite{BHJ}, for instance.
We then say that $w$ is {\em modulationally unstable}
if the $L^2_{2\pi}$-spectrum of $\mathcal{L}_\xi$ in the vicinity of the origin 
is not contained in the imaginary axis for some $|\xi|$ small, and it is modulationally stable otherwise.

The eigenvalue problem \eqref{E:Lxi} must in general be investigated numerically,
except in few special cases, e.g. completely integrable systems. 
In the case where $\lambda$ is near the origin and $|\xi|$ is small, nevertheless, 
we may take the perturbation theory approach in \cite{HJ2}, for instance, and address it analytically.
Specifically we first study the spectrum of the unmodulated operator $\mathcal{L}_0$
at the origin; note that zero is an eigenvalue of $\mathcal{L}_0$.
We then examine how the spectrum near the origin of the modulated operator $\mathcal{L}_\xi$
bifurcates from that of $\mathcal{L}_0$ for $|\xi|$ small.

In what follows, %$k>0$ is fixed and suppressed to simplify the exposition, unless specified otherwise.
thanks to Galilean invariance (see \eqref{E:galilean}), we may take $b=0$. 
%We write $\mathcal{L}_{\xi,a}$ for $\mathcal{L}_\xi(k,a,0)$. 

\

In the case of $a=0$, namely the zero solution, a straightforward calculation reveals that 
\[
\mathcal{L}_{\xi}(k,0,0) e^{inz}=i\omega_{n,\xi}e^{inz}
\quad \text{for all $n\in\mathbb{Z}$}\quad\text{for all $\xi\in[-1/2,1/2)$},
\]
where $\omega_{n,\xi}=(n+\xi)(m(k)-m(k(n+\xi)))$. 
%Consequently, $\text{spec}_{L^2_{2\pi}(\mathcal{L}_{\xi,0}) \subset i\mathbb{R}$ for all $\xi\in[-1/2,1/2)$ and hence the zero solution of \eqref{E:whitham} is spectrally stable to square integrable perturbations.
Notice that 
\[
\omega_{\pm1,0}=\omega_{0,0}=0.
\]
As a matter of fact, zero is an $L^2_{2\pi}$-eigenvalue of $\mathcal{L}_{0}(k,0,0)$ 
with algebraic and geometric multiplicity three, and
$\cos z$, $\sin z$ and $1$ form a basis of the corresponding eigenspace. 
For $|a|$ small, moreover zero is a generalized $L^2_{2\pi}$-eigenvalue of $\mathcal{L}_{0}(k,a,0)$
with algebraic multiplicity three and geometric multiplicity two, and 
\begin{align*}
\phi_1(z)=&\frac{1}{2(m(0)-m(k))}((\partial_bc)(\partial_aw)-(\partial_ac)(\partial_bw))(k,a,0)(z) \\
=&\cos z+\frac{-1/2+a\cos(2z)}{m(k)-m(2k)}+O(a^2), \\
\phi_2(z)=&-\frac1a(\partial_zw)(k,a,0)(z)=\sin z+a\frac{\sin(2z)}{m(k)-m(2k)}+O(a^2),\\
\phi_3(z)=&1
\end{align*}
form a basis of the corresponding generalized eigenspace. 
One may then proceed as in \cite{HJ2} by computing the $3\times 3$ matrix representing 
the action of $\mathcal{L}_{\xi,a}$ on the spectral space spanned by $\phi_j's$, $j=1,2,3$,
associated with three eigenvalues at zero. We summarize the conclusion below.
The detail is found in \cite{HJ2}, for instance, 
which is rather technical, involving commutators of $\M_k$ and $z$.
We merely pause to remark that an arbitrary $m(0)$ may be treated by a scaling argument.
In the case of $T/gd^2=1/3$, the quadratic order term in the Taylor expansion of $m$
vanishes, and the argument in \cite{HJ2} becomes inconclusive. 

%%%%%%%%%%%%%%%%%%%%%%%%%%%%%%%%%%%%%%%%%%%%%%%%%
%%%%%%%%%%%%%%%%%%%%%%%%%%%%%%%%%%%%%%%%%%%%%%%%%
%%%%%%%%%%%%%%%%%%%%%%%%%%%%%%%%%%%%%%%%%%%%%%%%%
\begin{theorem}[Modulational instability index]\label{t:mi}
Let $T>0$, $T/gd^2\neq 1/3$ and let $k\in\Sigma$.  
A $2\pi/k$-periodic traveling wave of \eqref{E:whitham}-\eqref{def:m}, \eqref{defT:m} 
with sufficiently small amplitude is modulationally unstable if $\Delta_{MI}(kd)<0$, where
\begin{equation}\label{def:index}
\Delta_{MI}(z)=\frac{(z m(z))''((z m(z))'-m(0))}{m(z)-m(2z)}\Delta_{BF}(z)
\end{equation}
and $\Delta_{BF}(z)=2(m(z)-m(2z))+((z m(z))'-m(0))$.
It is modulationally stable if $\Delta_{MI}(kd)>0$. 
\end{theorem}
%%%%%%%%%%%%%%%%%%%%%%%%%%%%%%%%%%%%%%%%%%%%%%%%%
%%%%%%%%%%%%%%%%%%%%%%%%%%%%%%%%%%%%%%%%%%%%%%%%%
%%%%%%%%%%%%%%%%%%%%%%%%%%%%%%%%%%%%%%%%%%%%%%%%%

Spectral instabilities away from the origin may be studied using a Krein signature method 
in \cite{DT}, for instance, although we do not consider them here.

\

Theorem~\ref{t:mi} identifies four mechanisms which cause change in the sign of 
the modulational instability index $\Delta_{MI}$, and hence
change in modulational stability and instability:  
\begin{itemize}
\item[(1)] $(z m(z))''=0$ at some $z=kd$, 
i.e. the group velocity $c_g(k)=((kd)m(kd))'$ attains an extremum at some wave number $k$; 
\item[(2)] $(z m(z))'=m(0)$ at some $z=kd$, 
i.e. the group velocidy coincides with the phase velocity $c_p(k)=m(kd)$ 
of the limiting long wave at $k=0$, resulting in the resonance between long and short waves;
\item[(3)] $m(z)=m(2z)$ at some $z=kd$, 
i.e. the phase velocities of the fundamental mode and the second harmonic coincide,
resulting in the ``second harmonic resonance";
note from the previous section that $k\in\Sigma$ excludes the possibility;
\item[(4)] $\Delta_{BF}(z)=0$ at some $z=kd$, 
which in the absence of the influence of surface tension occurs at $kd=1.146\dots$; see \cite{HJ2}.
\end{itemize}

Possibilities (1)-(3) are {\em purely dispersive} and 
they do not depending upon the specific nonlinearity of the governing equation.
As a matter of fact, they agree with those in \cite{kawahara} and \cite{DR}, for instance,
via formal multiple-scale expansions to the capillary-gravity wave problem. 
Furthermore they appear in the modulational instability indices for other, nonlinear dispersive equations; 
see \cite{HP1}, for instance.
The possibility (4), on the other hand, is associated with a resonance of dispersion and the nonlinearity,
and it depends upon the nonlinearity of the equation; see \cite{HP1}, for examples.
It does not agree with those in \cite{kawahara} and \cite{DR}.
This is not surprising since \eqref{E:whitham} does not include
many higher order nonlinearities of the physical problem.
It is interesting to attempt to utilize a more physically realistic nonlinearity,
but the task is well beyond the scope and interest of the present work. 

\

\begin{figure}[h]
\includegraphics[scale=0.7]{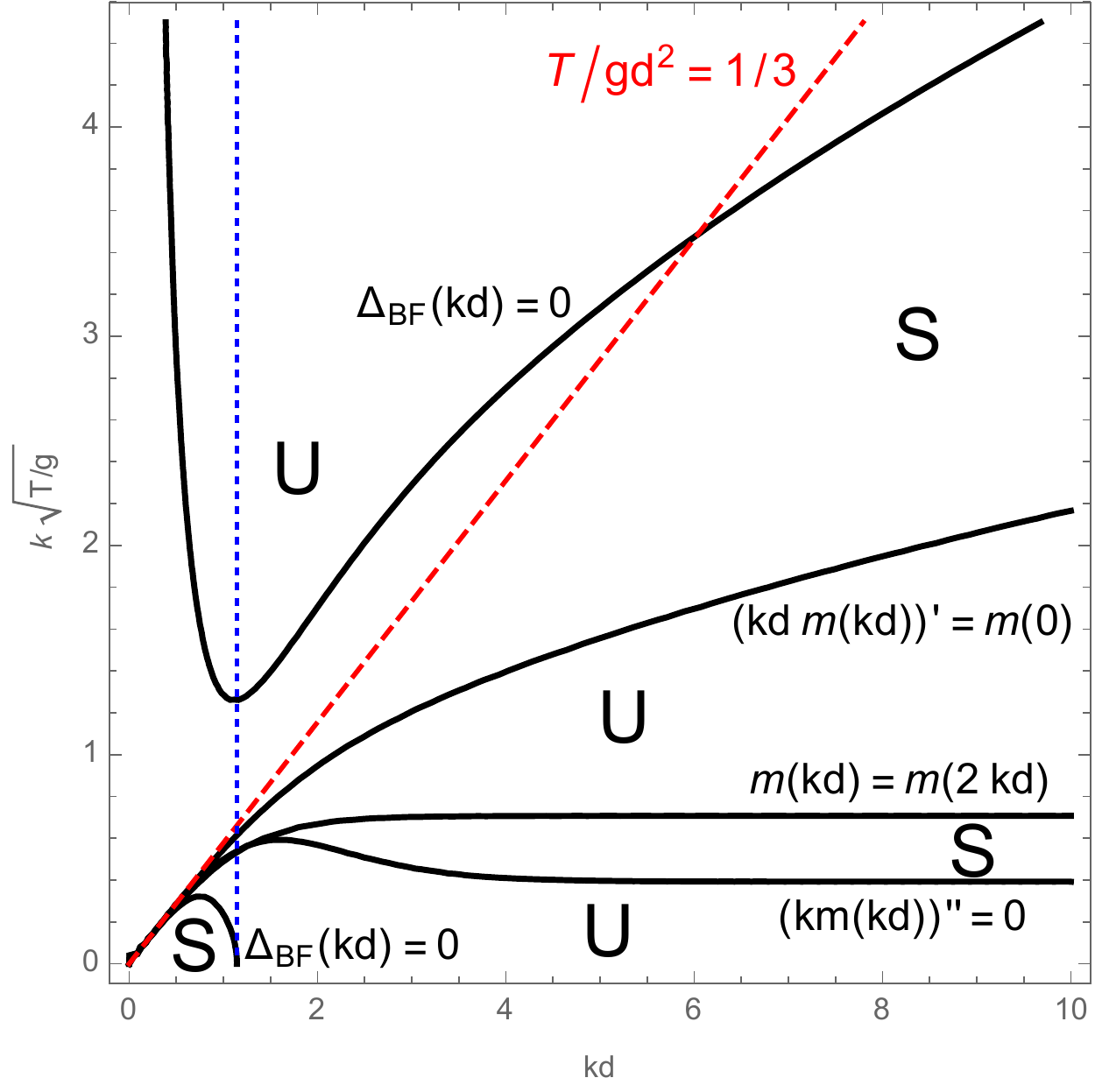}
\caption{Stability diagram in the $kd$ versus $k\sqrt{T/g}$ plane 
for small-amplitude, periodic traveling waves of the capillarity-modified Whitham equation. 
For $T\geq 0$ fixed, to interpret, one must envision a line through the origin with slope $\sqrt{T/gd^2}$.  
``S" and ``U" denote regions of modulational stability and instability.  
The vertical dashed line represents the positive root of $\Delta_{BF}$, 
in the absence of the effects of surface tension, approximately $1.146$, 
and the dashed line through the origin corresponds to the critical value $T/gd^2=1/3$. 
Other curves in the plane represent roots of $\Delta_{MI}$ and are labeled according to their mechanism.}\label{f:contour}
\end{figure}

The effects of surface tension considerably alter the modulational stability and instability 
in the capillarity-modified Whitham equation. Possibilities (1) through (4) may occur depending upon 
the coefficient of surface tension, the fluid depth and the carrier wave number,
and the sign of the modulational instability index is determined 
through a rather complicated balance of dispersion and the nonlinearity.
Figure~\ref{f:contour} illustrates regions in the $kd$ versus $k\sqrt{T/g}$ plane 
of modulational stability and instability.
It agrees with those in \cite{kawahara} and \cite{DR}, 
except the ``upper" branch of the curve along which $\Delta_{BF}=0$. 

\

%%%%%%%%%%%%%%%%%%%%%%%%%%%%%%%%%%%%%%%%%%%%%%%%%
%%%%%%%%%%%%%%%%%%%%%%%%%%%%%%%%%%%%%%%%%%%%%%%%%
%%%%%%%%%%%%%%%%%%%%%%%%%%%%%%%%%%%%%%%%%%%%%%%%%
%%%%%%%%%%%%%%%%%%%%%%%%%%%%%%%%%%%%%%%%%%%%%%%%%
%%%%%%%%%%%%%%%%%%%%%%%%%%%%%%%%%%%%%%%%%%%%%%%%%

In the gravity wave setting, i.e. $T=0$, a straightforward calculation reveals that
the phase velocity $m(z)$ and the group velocity $(zm(z))'$ monotonically decrease for all $z>0$ 
with no critical points, and hence possibilities (1) through (3) do not occur. 
Numerical evaluation, on the other hand, indicates that $\Delta_{BF}(z)$ makes one sign change, 
from positive to negative, at $z\approx1.145\dots$
and hence the possibility (4) contributes to change in the modulational stability and instability
; see \cite{HJ2} for the detail. 
Note in passing that it seems difficult to analytically study the sign of $\Delta_{BF}$, 
but it is explicit, involving hyperbolic functions, and hence amenable of numerical evaluation. 

When the coefficient of surface tension is large, or equivalently, the fluid depth is small so that $T/gd^2>1/3$, 
similarly, the phase velocity and the group velocity monotonically increase over the interval $(0,\infty)$ 
with no critical points, and hence possibilities (1) through (3) do not occur. 
Numerical evaluation indicates that $\Delta_{BF}$ makes one sign change,
and hence the possibility (4) contributes to change in the modulational stability and instability.
As a matter of fact, a small-amplitude, periodic traveling wave 
of \eqref{E:whitham}-\eqref{def:m}, \eqref{defT:m} is unstable to long wavelength perturbations, 
provided that the wave number is greater than a critical value,
and stable otherwise, qualitatively the same as in the gravity wave setting $T=0$. 

In the case of $0<T/gd^2<1/3$, in stark contrast,
the group velocity possesses a unique minimum. Moreover $(zm(z))'=m(0)=\sqrt{gd}$
and $m(z)=m(2z)$ each takes a root for some $z\in (0,\infty)$. Therefore 
each of possibilities (1) through (4) contributes to change in the modulational stability and instability, 
and we find intervals of stable and unstable wave numbers.
Note in particular that some short waves are modulationally stable.
Note however that for each $0<T/gd^2<1/3$, small-amplitude, periodic traveling waves 
become modulationally unstable as the wave number increases to infinity. 

\

In the capillary wave setting, i.e. $g=0$, we deduce from Figure~\ref{f:contour} that
small-amplitude, periodic traveling waves of \eqref{E:whitham}-\eqref{def:m}, \eqref{defT:m}
are modulationally unstable for all wave numbers. 
The result is in agreement with \cite{kawahara} 
but it goes against a conjecture by Barakat (see \cite{Barakat}, for instance)
that capillary waves in very shallow water are unlikely to be unstable. 

In the deep water limit as $kd\to\infty$, note in Figure~\ref{f:contour} that 
the curves along which $m(z)=m(2z)$ and $(zm(z))''=0$ approach the horizontal asymptotes 
at $1/\sqrt{2}$ and $\sqrt{2/\sqrt{3}-1}$, respectively, 
whereas the curve along which $(zm(z))'=m(0)$ and 
the ``upper" branch of the curve along which $\Delta_{BF}(z)=0$ grow without bounds. 
Therefore possibilities (1) and (3) contribute to change in the modulational stability and instability. 
Specifically, a small-amplitude capillary-gravity wave in deep water 
is expected modulationally stable in the range $2/\sqrt{3}-1<Tk^2/g<1/2$, and unstable otherwise.
The result agrees with that in \cite{Lighthill-deep}, for instance, via formal argument for the physical problem.
The limiting dispersion symbol $\sqrt{T|k|+g/|k|}$ is not smooth, however, 
and the present treatment, or that in \cite{HJ2}, fails. 

\

In neighborhoods of boundaries of stable and unstable regions in Figure~\ref{f:contour}, 
the present treatment fails and one must either go to higher order or consider a different method. 
The resonance between long and short waves was studied in \cite{DR}, for instance. 
The second harmonic resonance has been studied by many researchers, albeit formally; 
see, for instance, the book by Craik~\cite{Craik}.

%%%%%%%%%%%%%%%%%%%%%%%%%%%%%%%%%%%%%%%%%%%%%%%%%
%%%%%%%%%%%%%%%%%%%%%%%%%%%%%%%%%%%%%%%%%%%%%%%%%
%%%%%%%%%%%%%%%%%%%%%%%%%%%%%%%%%%%%%%%%%%%%%%%%%
\section{Effects of constant vorticity}\label{sec:vorticity}
%%%%%%%%%%%%%%%%%%%%%%%%%%%%%%%%%%%%%%%%%%%%%%%%%
%%%%%%%%%%%%%%%%%%%%%%%%%%%%%%%%%%%%%%%%%%%%%%%%%
%%%%%%%%%%%%%%%%%%%%%%%%%%%%%%%%%%%%%%%%%%%%%%%%%

We extend the results of previous sections 
to the {\em vorticity-modified Whitham equation}, \eqref{E:whitham}-\eqref{def:m}, where, abusing notation,
\begin{equation}\label{defw:m}
m_{\gamma, \pm}(k)=\gamma d\frac{\tanh(k)}{2k}
+\sqrt{gd\frac{\tanh(k)}{k}+\gamma^2d^2\frac{\tanh^2(k)}{4k^2}}.
\end{equation}
As a matter of fact, $kd\mapsto k$ turns \eqref{def:mw} to \eqref{defw:m}.
Note that $\pm m_{\gamma,\pm}(k)>0$ for all $\gamma\in\mathbb{R}$ for all $k>0$. 
Note moreover that $m_{\gamma,\pm}(k)>m_{0,\pm}(k)$ for all $k>0$ if $\gamma>0$ and 
$m_{\gamma,\pm}(k)<m_{0,\pm}(k)$ for all $k>0$ if $\gamma<0$. 
In other words, a positive constant vorticity increases the phase velocity 
and a negative constant vorticity decreases. %; see also the discussion in \cite{CKS15}. 

One may assume, without loss of generality, the $+$ sign in \eqref{defw:m}.
As a matter of fact, 
\begin{equation}\label{E:vorsymmetry}
m_{\gamma,-}(k)=m_{-\gamma,+}(k) \qquad \text{for all $\gamma\in\mathbb{R}$ for all $k>0$.}
\end{equation}
In what follows, we shall take the $+$ sign in \eqref{defw:m} and suppress it to simplify the exposition, 
unless specified otherwise.
In the long wave limit as $k\to 0$, a KdV equation whose coefficients depend upon $\gamma$
approximates up to second order the dispersion relation of the vorticity-modified Whitham equation,
and hence the water wave problem with constant vorticity. 
%\[
%u_t+\Big(\Big(\frac{\gamma d}{2}+\sqrt{\frac{\gamma^2d^2}{4}+gd}\Big)
%-\frac16\Big(\frac{\gamma d}{2}+\frac{\frac12\gamma^2d^2+gd}{\sqrt{\gamma^2d^2+gd}}
%\partial_x^2\Big)\Big)u+uu_x=0.
%\]
In the short wave limit as $k\to \infty$, on the other hand, \eqref{E:whitham}-\eqref{def:m} and \eqref{defw:m}
becomes, after normalization of parameters, a perturbation of 
\[
u_t+\Lambda^{-1/2}u_x+uu_x=0,\qquad\text{where}\quad 
\widehat{\Lambda f}(k)=|k|\widehat{f}(k).
\]
As a matter of fact, it was argued in \cite{Hur-breaking} to have relevances to 
the water wave problem in two dimensions in the infinite depths. 
In particular, it shares in common with the physical problem the dispersion relation and scaling symmetry.

\

\begin{figure}[h]
\includegraphics[scale=0.5]{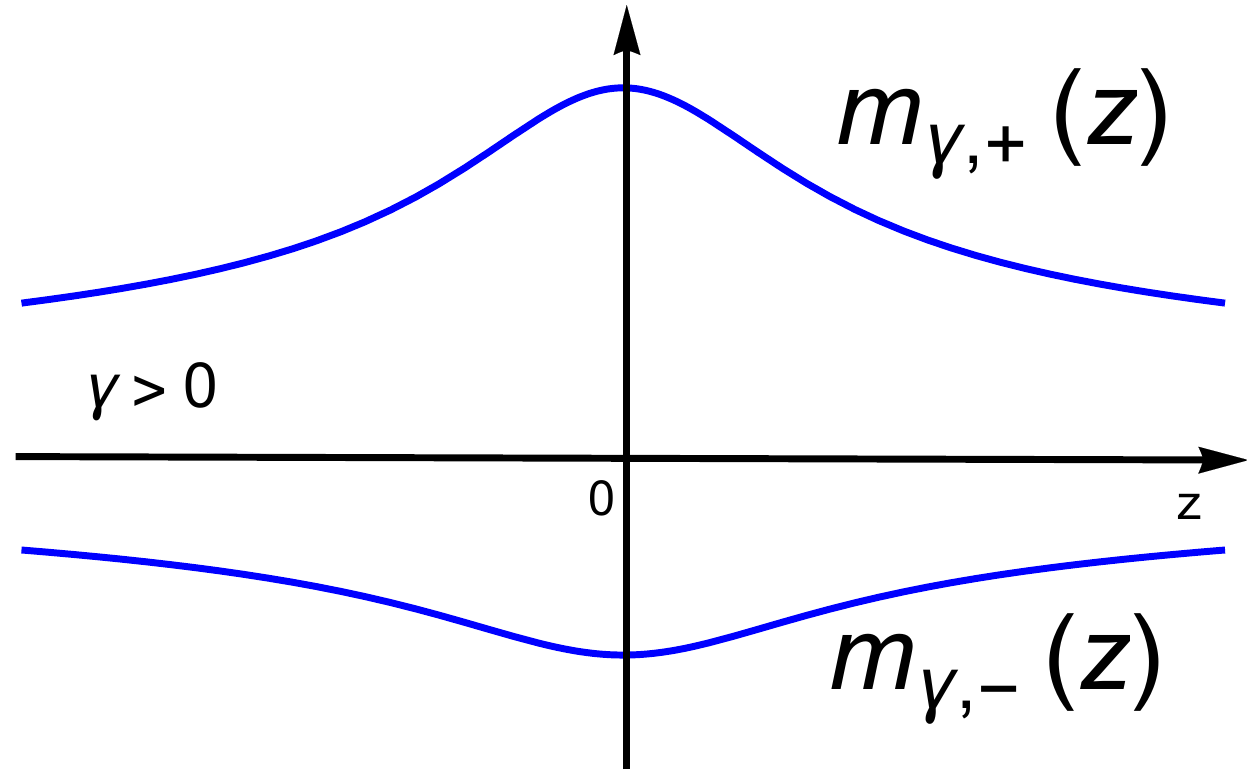}
\caption{Schematic plot of the symbols $m_{\gamma,\pm}$ for a positive constant vorticity.}\label{f:vorsymbol}
\end{figure}

For each $\gamma\in\mathbb{R}$, 
note that $m(z)$ is smooth, even and monotonically decreasing to zero as $z\to \infty$;
see Figure~\ref{f:vorsymbol}. Therefore $m(k)>m(kN)$ for all $k>0$ for all $N=2,3,\dots$. 
We then follow the arguments in Section~\ref{sec:existence}, or in \cite{HJ2}, to obtain 
for each $k>0$ a smooth, two-parameter family of small-amplitude, $2\pi/k$-periodic waves
$w=w(k,a,b)$ of \eqref{E:whitham}-\eqref{def:m} and \eqref{defw:m}
traveling at the speed $c=c(k,a,b)$, where $|a|$ and $|b|$ are sufficiently small.
Moreover $w$ and $c$ satisfy, respectively, \eqref{E:w} and \eqref{E:c}.

%%%%%%%%%%%%%%%%%%%%%%%%%%%%%%%%%%%%%%%%%%%%%%%%%
%%%%%%%%%%%%%%%%%%%%%%%%%%%%%%%%%%%%%%%%%%%%%%%%%
%%%%%%%%%%%%%%%%%%%%%%%%%%%%%%%%%%%%%%%%%%%%%%%%%
%%%%%%%%%%%%%%%%%%%%%%%%%%%%%%%%%%%%%%%%%%%%%%%%%
 
 \
 
Furthermore we follow the arguments in Section~\ref{sec:MI}, or in \cite{HJ2},
to derive the modulational instability index $\Delta_{MI}$ in \eqref{def:index},
which determines the modulational stability and instability  
of a small-amplitude, periodic traveling wave of \eqref{E:whitham}-\eqref{def:m} and \eqref{defw:m}.

For any $\gamma\in \mathbb{R}$, a straightforward calculation reveals that
the phase velocity $m(z)$ and the group velocity $(zm(z))'$ monotonically decrease for all $z>0$
without critical points, and hence possibilities (1) through (3) following Theorem~\ref{t:mi} do not occur.
Numerical evaluation, on the other hand, indicates that 
$\Delta_{BF}(z)$ in \eqref{def:index} makes one sign change, 
from positive to negative, at $z_c(\gamma)$, say, and hence the possibility (4) 
following Theorem~\ref{t:mi} contributes to change in the modulational stability and instability.
Therefore a small-amplitude, $2\pi/k$-periodic traveling wave of
\eqref{E:whitham}-\eqref{def:m} and \eqref{defw:m} is unstable to long wavelength perturbations 
if $kd>z_c$, and stable otherwise, 
qualitatively the same as in the zero vorticity setting $\gamma=~0$. 
Furthermore, a straightforward calculation reveals that 
$z_c(\gamma)$ decreases with $\gamma$ and 
\[
\lim_{\gamma\to-\infty}z_{c}(\gamma)=\infty\quad\text{and}\quad
\lim_{\gamma\to\infty}z_{c}(\gamma)=0.957\ldots.
\]

In the case of the $-$ sign in \eqref{defw:m}, thanks to \eqref{E:vorsymmetry}, we obtain the same result.
As a matter of fact, 
\[
\Delta_{MI,+}(z;\gamma)=\Delta_{MI,-}(z,\gamma)
\] 
for all $z>0$ for all $\gamma\in\mathbb{R}$. But
\[
\Delta_{\rm BF,-}(z;\gamma)=-\Delta_{\rm BF,+}(z;-\gamma)
\]
for all $z>0$ for all $\gamma\in\mathbb{R}$, and therefore 
the critical wave number in the case of the $-$ sign in \eqref{defw:m} increases
with $\gamma$ from $0.957\dots$ to infinity. 

\begin{figure}[h]
\includegraphics[scale=0.5]{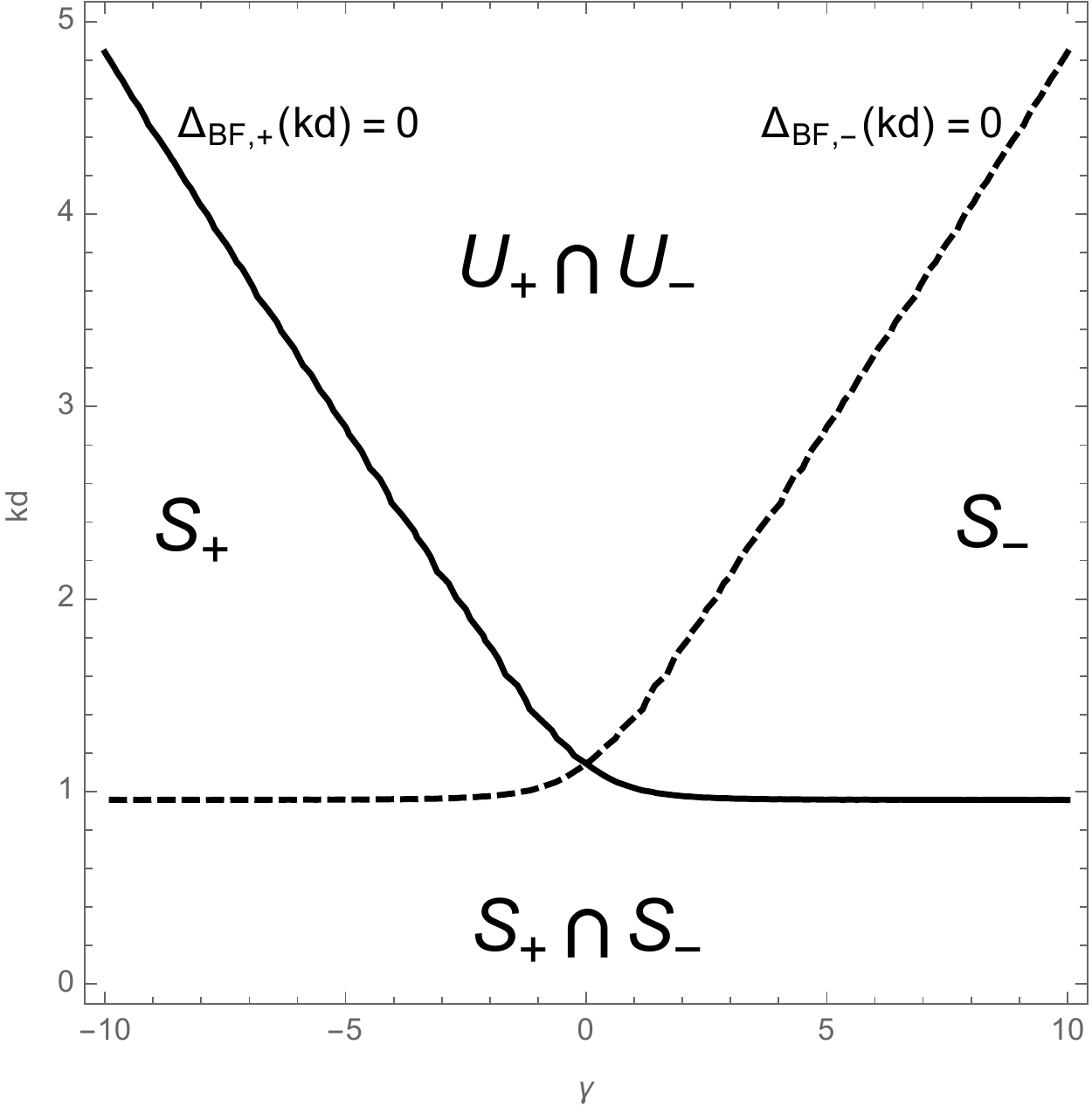}
\caption{
Stability diagram in the $\gamma$ versus $kd$ plane 
for small-amplitude, periodic traveling wave of the vorticity-modified Whitham equation.
The solid and dashed curves correspond to the level curves 
$\Delta_{BF,+}(z)=0$ and $\Delta_{BF,-}(z)=0$, respectively.   
$``S_+"$ and $``S_-"$ denote regions of modulational stability 
with the corresponding choice of the sign in \eqref{defw:m}, 
and $``U_+\cap U_-"$ is the intersection of regions of instability for two equations. 
Two level sets intersect at $\gamma=0$ and $kd=1.145\dots$.}\label{f:vorcontour}
\end{figure}

To recapitulate, for each $\gamma \in \mathbb{R}$, 
a small-amplitude, periodic traveling wave of \eqref{E:whitham}-\eqref{def:m} and \eqref{defw:m}
is modulationally unstable, provided that the wave number is greater than a critical number, 
and stable otherwise, similarly to the zero vorticity setting. 
The critical wave number monotonically increases to infinity
as $\gamma\to-\infty$ in the case of the $+$ sign in \eqref{defw:m} and 
as $\gamma\to+\infty$ in the case of the $-$ sign in \eqref{defw:m},
and it is bounded below by $0.957\dots$.
Figure~\ref{f:vorcontour} illustrates regions in the $\gamma$ versus $kd$ plane
of modulational stability and instability.

%%%%%%%%%%%%%%%%%%%%%%%%%%%%%%%%%%%%%%%%%%%%%%%%%
%%%%%%%%%%%%%%%%%%%%%%%%%%%%%%%%%%%%%%%%%%%%%%%%%
%%%%%%%%%%%%%%%%%%%%%%%%%%%%%%%%%%%%%%%%%%%%%%%%%
%%%%%%%%%%%%%%%%%%%%%%%%%%%%%%%%%%%%%%%%%%%%%%%%%
%%%%%%%%%%%%%%%%%%%%%%%%%%%%%%%%%%%%%%%%%%%%%%%%%
%%%%%%%%%%%%%%%%%%%%%%%%%%%%%%%%%%%%%%%%%%%%%%%%%

In particular, a small-amplitude, periodic traveling wave of the vorticity-modified Whitham equation
can be made modulationally stable or unstable, by ensuring that $|\gamma|$ is sufficiently large.  
For example, consider some $kd>1.145\ldots$, so that 
the corresponding, small-amplitude periodic traveling wave 
of the Whitham equation (in the zero vorticity setting) is modulationally unstable. 
It can be transitioned into the region of modulational stability 
by taking $\gamma$ sufficiently large and negative in the case of the $+$ sign in \eqref{defw:m}
and positive in the case of the $-$ sign in \eqref{defw:m}.
Similarly, $kd$ in the range $(0.957\ldots,1.145\ldots)$, so that the corresponding 
small-amplitude, periodic traveling wave of the Whitham equation is modulationally stable,
can be transitioned to the region of modulational instability
by taking $\gamma$ sufficiently large and positive in the case of $+$ sign in \eqref{defw:m}
and negative in the case of $-$ sign in \eqref{defw:m}.
Note that all small-amplitude, $2\pi/k$-periodic traveling waves of \eqref{E:whitham}-\eqref{def:m}
and \eqref{defw:m} are modulationally stable if $0<kd<0.957\dots$, 
for all $\gamma\in\mathbb{R}$ regardlessly of the sign in \eqref{defw:m}.

Recent numerical experiments in \cite{OSV} suggest that 
a small-amplitude, periodic traveling wave in water with constant vorticity
can be made modulationally stable by ensuring that the constant vorticity is sufficiently large in strength.
Our results bear it out, observing modulational stabilization 
in each branch of the phase velocity in the vorticity-modified Whitham equation
for $\gamma$ sufficiently large and with an appropriate sign (depending upon the branch). 

Taking the regions of modulational stability in Figure~\ref{f:vorcontour} together suggests that 
if one were to consider \emph{both branches} of the phase velocity simultaneously, say, 
in a bi-directional Whitham equation, then
one would be likely to observe the modulational stabilization effects 
of constant vorticity, independent of the sign of $\gamma$, as reported in \cite{OSV}.
This is an interesting direction for future investigations.

\

In the deep water limit as $kd\to\infty$, Figure~\ref{f:vorcontour} indicates that 
small-amplitude, periodic traveling waves of \eqref{E:whitham}-\eqref{def:m} and \eqref{defw:m}
be modulationally unstable for all wave numbers for all vorticities. 
In the zero vorticity setting, this is in agreement with the results in \cite{BF} and \cite{Whitham1967}.
The limiting symbol $\sqrt{g/|k|}$ is not smooth, however, 
and the present treatment, or that in \cite{HJ2}, fails.
 
%%%%%%%%%%%%%%%%%%%%%%%%%%%%%%%%%%%%%%%%%%%%%%%%%
%%%%%%%%%%%%%%%%%%%%%%%%%%%%%%%%%%%%%%%%%%%%%%%%%
\subsection*{Acknowledgements}
VMH is supported by the National Science Foundation grant CAREER DMS-1352597, 
an Alfred P. Sloan Foundation fellowship, and a Beckman fellowship 
of the Center for Advanced Study at the University of Illinois at Urbana-Champaign. 
MAJ is supported by the National Science Foundation under grant DMS-1211183.
The authors thank the anonymous referees for their careful reading of the manuscript
and many helpful comments and references.	 

\bibliographystyle{amsalpha}
\bibliography{stabilityBib}

\providecommand{\bysame}{\leavevmode\hbox to3em{\hrulefill}\thinspace}
\providecommand{\MR}{\relax\ifhmode\unskip\space\fi MR }
\providecommand{\MRhref}[2]{%
  \href{http://www.ams.org/mathscinet-getitem?mr=#1}{#2}
}
\providecommand{\href}[2]{#2}
\begin{thebibliography}{SKCK14}

\bibitem[Air45]{Airy}
G.~B. Airy, \emph{Tides and waves}, Encyclopedia Metropolitana \textbf{5}
  (1845), 241--396.

\bibitem[Bar]{Barakat}
Richard Barakat, \emph{Instability of periodic capillary-gravity waves on deep
  water}, Wave Motion \textbf{6}, no.~2, 155--165.

\bibitem[BF67]{BF}
T.~B. Benjamin and J.~E. Feir, \emph{The disintegration of wave trains on deep
  water. {P}art 1. {T}heory}, J. Fluid Mech. \textbf{27} (1967), no.~3,
  417--437.

\bibitem[BHJ16]{BHJ}
Jared~C. Bronski, Vera~Mikyoung Hur, and Mathew~A. Johnson, \emph{Modulational
  instability in equations of {K}d{V} type}, New Approaches to Nonlinear Waves,
  Lecture Notes in Physics, Springer, 2016.

\bibitem[BKN13]{BKN}
Handan Borluk, Henrik Kalisch, and David~P. Nicholls, \emph{A numerical study
  of the {W}hitham equation as a model for steady surface water waves},
  preprint (2013).

\bibitem[BN67]{BN}
D.~J. Benney and A.~C. Newell, \emph{The propagation of nonlinear wave
  envelopes}, J. Math. and Phys. \textbf{46} (1967), 133--139. \MR{0241052 (39
  \#2397)}

\bibitem[CDG95]{CDG95}
T.~Colin, F.~Dias, and J.-M. Ghidaglia, \emph{On rotational effects in the
  modulations of weakly nonlinear water waves over finite depth}, European J.
  Mech. B Fluids \textbf{14} (1995), no.~6, 775--793. \MR{1364731 (97i:76019)}

\bibitem[CE98]{CE-breaking}
Adrian Constantin and Joachim Escher, \emph{Wave breaking for nonlinear
  nonlocal shallow water equations}, Acta Math. \textbf{181} (1998), no.~2,
  229--243. \MR{1668586 (2000b:35206)}

\bibitem[Cra85]{Craik}
Alex D.~D. Craik, \emph{Wave interactions and fluid flows}, Cambridge
  Monographs on Mechanics and Applied Mathematics, Cambridge University Press,
  Cambridge, 1985. \MR{896268 (88f:76001)}

\bibitem[CS04]{CS2004}
Adrian Constantin and Walter Strauss, \emph{Exact steady periodic water waves
  with vorticity}, Comm. Pure Appl. Math. \textbf{57} (2004), no.~4, 481--527.
  \MR{2027299 (2004i:76018)}

\bibitem[DR77]{DR}
V.~D. Djordjevi{\'c} and L.~G. Redekopp, \emph{On two-dimensional packets of
  capillary-gravity waves}, J. Fluid Mech. \textbf{79} (1977), no.~4, 703--714.
  \MR{0443555 (56 \#1924)}

\bibitem[DT15]{DT}
Bernard Deconinck and Olga Trichtchenko, \emph{High-frequency instabilities of
  small-amplitude solutions of hamiltonian pdes}, Preprint (2015).

\bibitem[EEP15]{EEP}
Mats Ehrnstr{\"o}m, Joachim Escher, and Long Pei, \emph{A note on the local
  well-posedness for the {W}hitham equation}, Elliptic and Parabolic Equations,
  Springer Proceedings in Mathematics \& Statistics, vol. 119, Springer
  International Publishing, 2015, pp.~63--75.

\bibitem[EGW12]{EGW}
Mats Ehrnstr{\"o}m, Mark~D. Groves, and Erik Wahl{\'e}n, \emph{On the existence
  and stability of solitary-wave solutions to a class of evolution equations of
  {W}hitham type}, Nonlinearity \textbf{25} (2012), no.~10, 2903--2936.
  \MR{2979975}

\bibitem[Ehr15]{Mats}
Mats Ehrnstr\"om, \emph{On {W}hitham's conjecture of a highest cusped wave for
  a nonlocal shallow water wave equation}, Mathematisches Forschungsinstitut
  Oberwolfach Report No. 19/2015, 2015.

\bibitem[EK09]{EK}
Mats Ehrnstr{\"o}m and Henrik Kalisch, \emph{Traveling waves for the {W}hitham
  equation}, Differential Integral Equations \textbf{22} (2009), no.~11-12,
  1193--1210. \MR{2555644 (2010k:35403)}

\bibitem[EK13]{EK2}
M.~Ehrnstr{\"o}m and H.~Kalisch, \emph{Global bifurcation for the {W}hitham
  equation}, Math. Model. Nat. Phenom. \textbf{8} (2013), no.~5, 13--30.
  \MR{3123360}

\bibitem[HJ15]{HJ2}
Vera~Mikyoung Hur and Mathew~A. Johnson, \emph{Modulational instability in the
  {W}hitham equation for water waves}, Stud. Appl. Math. \textbf{134} (2015),
  no.~1, 120--143.

\bibitem[HL08]{HL08}
Vera~Mikyoung Hur and Zhiwu Lin, \emph{Unstable surface waves in running
  water}, Comm. Math. Phys. \textbf{282} (2008), no.~3, 733--796. \MR{2426143
  (2009d:76023)}

\bibitem[HO72]{HO}
Hidenori Hasimoto and Hiroaki Ono, \emph{Nonlinear modulation of gravity
  waves}, J. Phys. Soc. Japan \textbf{33} (1972), no.~3, 805--811.

\bibitem[Hog86]{Hogan}
S.~J. Hogan, \emph{The potential form of the fourth-order evolution equation
  for deep-water gravity-capillary waves}, Phys. Fluids \textbf{29} (1986),
  no.~10, 3479--3480. \MR{861945 (87i:76013)}

\bibitem[HP15]{HP1}
Vera~Mikyoung Hur and Ashish~Kumar Pandey, \emph{Modulational instability in
  nonlinear dispersive equations}, in preparation (2015).

\bibitem[Hur12]{Hur-breaking}
Vera~Mikyoung Hur, \emph{On the formation of singularities for surface water
  waves}, Commun. Pure Appl. Anal. \textbf{11} (2012), no.~4, 1465--1474.
  \MR{2900797}

\bibitem[Hur15]{Hur-breaking1}
Vera~Mikyoung Hur, \emph{Breaking in the {W}hitham equation for shallow water
  waves}, arxiv:1506.04075 (2015).

\bibitem[Jon89]{Jones}
M.~C.~W. Jones, \emph{Small amplitude capillary-gravity waves in a channel of
  finite depth}, Glasgow Math. J. \textbf{31} (1989), no.~2, 141--160.
  \MR{997809 (90h:76025)}

\bibitem[Kar12]{Kar12}
Paschalis Karageorgis, \emph{Dispersion relation for water waves with
  non-constant vorticity}, Eur. J. Mech. B Fluids \textbf{34} (2012), 7--12.
  \MR{2927940}

\bibitem[Kat76]{K}
Tosio Kato, \emph{Perturbation theory for linear operators}, second ed.,
  Springer-Verlag, Berlin-New York, 1976, Grundlehren der Mathematischen
  Wissenschaften, Band 132. \MR{0407617 (53 \#11389)}

\bibitem[Kaw75]{kawahara}
Takuji Kawahara, \emph{Nonlinear self-modulation of capillary-gravity waves on
  liquid layer}, J. Phys. Soc. Japan \textbf{38} (1975), no.~1, 265--270.
  \MR{678043 (83k:76081)}

\bibitem[Lig65]{Lighthill-deep}
M.~J. Lighthill, \emph{Group velocity}, J. Inst. Math. Appl. \textbf{1} (1965),
  1--28. \MR{0184458 (32 \#1930)}

\bibitem[Mar13]{M2013}
Calin~Iulian Martin, \emph{Local bifurcation and regularity for steady periodic
  capillary-gravity water waves with constant vorticity}, Nonlinear Analysis:
  Real World Applications (2013), 131--149.

\bibitem[MKD14]{MKD}
Daulet Moldabayev, Henrik Kalisch, and Denys Dutykh, \emph{The {W}hitham
  equation as a model for surface water waves}, preprint (2014).

\bibitem[NS94]{NS}
P.~I. Naumkin and I.~A. Shishmar{\"e}v, \emph{Nonlinear nonlocal equations in
  the theory of waves}, Translations of Mathematical Monographs, vol. 133,
  American Mathematical Society, Providence, RI, 1994, Translated from the
  Russian manuscript by Boris Gommerstadt. \MR{1261868 (94m:35230)}

\bibitem[OSV15]{OSV}
K.~Oliveras, P.~Sprenger, and V.~Vasan, \emph{Stability of periodic traveling
  waves with constant vorticity}, preprint (2015).

\bibitem[Sim69]{Simmons}
William~Frederick Simmons, \emph{A variational method for weak resonant wave
  interactions}, Proc. R. Soc. Lond. Ser. A Math. Phys. Sci. \textbf{309}
  (1969), no.~1499, 551--557.

\bibitem[SKCK14]{SKCK2014}
Nathan Sanford, Keri Kodama, John~D. Carter, and Henrik Kalisch,
  \emph{Stability of traveling wave solutions to the {W}hitham equation}, Phys.
  Lett. A \textbf{378} (2014), no.~30-31, 2100--2107. \MR{3226084}

\bibitem[TdSP88]{SP}
A.~F. Teles~da Silva and D.~H. Peregrine, \emph{Steep, steady surface waves on
  water of finite depth with constant vorticity}, J. Fluid Mech. \textbf{195}
  (1988), 281--302. \MR{985439 (90a:76061)}

\bibitem[Whi67]{Whitham1967}
G.~B. Whitham, \emph{Non-linear dispersion of water waves}, J. Fluid Mech.
  \textbf{27} (1967), 399--412. \MR{0208903 (34 \#8711)}

\bibitem[Whi74]{Whitham}
\bysame, \emph{Linear and nonlinear waves}, Pure and Applied Mathematics (New
  York), Wiley-Interscience [John Wiley \& Sons], New York, 1974. \MR{0483954
  (58 \#3905)}

\bibitem[Zak68]{Zakharov-WW}
V.~E. Zakharov, \emph{Stability of periodic waves of finite amplitude on the
  surface of a deep fluid}, J. Appl. Mech. Tech. Phys. \textbf{9} (1968),
  no.~2, 190--194.

\end{thebibliography}

\end{document}